\documentclass{amsart}

\usepackage{amsmath,amsthm,amscd,amssymb,eucal,mathrsfs}
\usepackage{amsfonts}
\usepackage{latexsym}
\usepackage{graphicx}
\usepackage{multicol}

\usepackage[all]{xy}

\newtheorem{thm}{Theorem}[section]
\newtheorem{prop}[thm]{Proposition}
\newtheorem{lem}[thm]{Lemma}
\newtheorem{cor}[thm]{Corollary}
\newtheorem{rem}[thm]{Remark}
\newtheorem{dfn}[thm]{Definition}

\newtheorem{exa}[thm]{Example}

\newcommand{\mathset}[1]{{\{#1\}}} 

\newcommand{\norm}[1]{\left\|#1\right\|}

\DeclareMathOperator{\Val}{Val}

\newcommand{\Mumf}{{\rm Mumf}}
\DeclareMathOperator{\PGL}{PGL}

\DeclareMathOperator{\branch}{br}

\newcommand{\An}{{\rm an}}
\newcommand{\trop}{{\rm trop}}
\DeclareMathOperator{\Div}{Div}
\DeclareMathOperator{\Spec}{Spec}
\newcommand{\closure}{{\rm cl}}

\begin{document}

\title{$p$-adic Hurwitz numbers}

\author{Patrick Erik Bradley}
 
\date{\today}

\begin{abstract}
We introduce stable tropical curves, and use these to count covers of
the $p$-adic projective line of fixed degree and ramification types by
Mumford curves in terms of tropical Hurwitz numbers. Our counts depend
on the branch loci of the covers. 
\end{abstract}

\maketitle

\section{Introduction}

The recent success of tropical geometry lies in its ability to provide
elementary methods in the enumerative geometry of algebraic curves.
For example, plane Gromov-Witten invariants can be defined in the tropical
setting in order to obtain the classical recursion formulae \cite{GathmannMarkwigWDVV,GathmannMarkwigCH}. 

By its own nature, it seems straightforward to relate tropical geometry to 
non-archimedean geometry and obtain an interpretation of e.g.\ the 
$j$-invariant of elliptic curves as the cycle length of  their tropicalisations
\cite{KMM}. This interpretation becomes strongly supported by the
tropicalisation of $p$-adic analytic spaces in \cite{Gubler1998}.
The latter, in fact, has its precedent in the use of tropical tori
in order to study the Jacobian of $p$-adic Mumford curves
in \cite{GvP1980}.

Important objects in the enumerative geometry of curves are moduli spaces
and their compactifications. These spaces exist also in tropical geometry,
although they are defined mostly in an ad-hoc fashion and not always are compact. Hence, a foundational approach to moduli spaces via something like 
``tropical stacks'' would be desireable. A first step in this
direction is found in the upcoming book \cite{Mikhalkin},
where tropical varieties are defined.
This should lead into proving theorems about tropicalised
spaces (or stacks)
by ``tropicalising'' existing
proofs for their classical or $p$-adic counterparts. 

The aim of this article is to relate Hurwitz numbers for covers
of the $p$-adic projective line by Mumford curves to tropical Hurwitz numbers.
The fact whether or not the upper curve is a Mumford curve
depends on the branch locus in a not so obvious way. First results
have been found for hyperelliptic covers in \cite{GvP1980}.
The results were extended to cyclic covers \cite{Ulrich1981,Brad-manmat124}
and to general finite Galois covers \cite{Brad-MathZ251}. In the latter
case, the motivation was inverse Galois theory, and Hurwitz
spaces for Galois covers were used. However, the situation
for covers which are not Galois remained unclear. It is 
the use of tropical geometry here which allows in principle to
count degree $d$ covers of $\mathbb{P}^1$ by genus $g$ Mumford curves
and fixed ramification types $\eta_0,\eta_1,\dots$,
depending on the branch locus. 
The most important ingredient
for actually deciding whether a given tropical cover comes from a Mumford
curve is a tropical version of the Riemann-Hurwitz formula 
(Theorem \ref{RHtrop}) which yields the characterising
 criterion that the ramification divisor be effective (Corollary \ref{MRH}).
Our main enumerative result is that the classical Hurwitz number
$H_d^g(\eta_0,\dots,\eta_{n+1})$ counts covers by Mumford curves,
if the branch tree is binary (Theorem \ref{bin=mumf}).
In particular,  the
double Hurwitz numbers studied by Cavalieri et al.\ \cite{CJM} 
count Mumford curves,
as their tropical covers correspond to tropical maps 
with effective ramification 
above comb-shaped binary trees
(Corollary \ref{comb=mumf}).

We do not refrain from using ad-hoc constructions
for the moduli spaces of stable tropical curves and maps.
However, a more foundational approach is the content of ongoing work.
We benefit
from the definition of projective tropical variety in \cite{Mikhalkin}
which is in fact a straightforward imitation of the corresponding
classical notion. In this sense, our $n$-pointed
tropical curves are tropically quasi-projective,
and the moduli spaces we introduce here are most likely projective.

It turns out that the philosophy of tropical curves being generically 
tropicalisations of Mumford curves 
can be proven also in the case of Hurwitz covers.

\medskip
The author acknowledges support from
the
 Deutsche Forschungsgemeinschaft through the
project BR 3513/1-1.

\section{Generalities}
Let $p$ be a prime number. It will hardly ever be referred to.
Our ground field will be $K=\mathbb{C}_p$, the  
completion of the algebraic closure of the field $\mathbb{Q}_p$
of $p$-adic numbers.  
The valuation map
will be denoted by $v\colon K\to \mathbb{R}\cup\mathset{\infty}$.

\section{Tropicalised stable  maps to $\mathbb{P}^1$}

Here, we introduce the notion of stable tropical curve and stable tropical
map to $\mathbb{P}^1$. These notions generalise the existing definitions
from \cite{GathmannMarkwigWDVV,CJM}.

\subsection{Tropicalisation of quasi-projective curves}

Let $\iota\colon C \to\mathbb{P}^N$ be the embedding of a $p$-adic 
quasi-projective curve into the projective space $\mathbb{P}^N$.
We will define a tropicalisation map for $\iota$ generalising the
approach of \cite{Payne}.

\medskip
There is the extended valuation map on affine space $\mathbb{A}^n$:
$$
v\colon\mathbb{A}^N\to \mathbb{TA}^N:=(\mathbb{TA}^1)^N,
\;
(x_1,\dots,x_N)\mapsto (v(x_1),\dots,v(x_N)),
$$ 
which has been used in \cite{Payne} to study Berkovich-analytifications
of affine curves.

The next step is to glue the maps $v_i=v\colon U_i\to \mathbb{TA}^N$
on the standard affine charts $U_i\cong\mathbb{A}^N$ of projective
$N$-space $\mathbb{P}^N$. This yields our tropicalisation map
$$
\Val\colon\mathbb{P}^N\to \mathbb{TP}
$$
to the tropical projective space from \cite{Mikhalkin}. The space
$\mathbb{TP}^N$ is obtained by glueing
the tropical affine pieces $\mathbb{TA}^N$ 
in the ``tropically'' analogous way as usual projective
space $\mathbb{P}^N$ is constructed from $\mathbb{A}^N$.

Let $(C,p_1,\dots,p_n)$ be a smooth irreducible $n$-punctured
projective curve of
genus $g$ defined over $K$ 
with $n$ $K$-rational punctures $p_1,\dots,p_n$.
Let $\iota\colon C\to\mathbb{P}^N$ be the pluricanonical
embedding of the curve given by $\omega_C(p_1+\dots+p_n)^{\otimes e}$,
 a sufficiently high power of
the ample line bundle $\omega_C(p_1+\dots+p_n)$. It has the
property 
\begin{align}
O(1)|_C=\omega_C(p_1+\dots+p_n)^{\otimes e}. \label{can=taut}
\end{align}

The closure $\Gamma:=\trop_\iota(C)$ of  $\Val(\iota(C))$ is known to be
 an abstract tropical curve 
with $n$ punctures and first Betti number at most $g$, and the punctures
are represented by unbounded  edges of $\Gamma$
(cf.\ the following subsection). From the point of view
of Berkovich \cite[\S 4]{Berkovich1990}, the tropical curve $\Gamma$
coincides with the skeleton of $C^\An$, and the unbounded edges of $\Gamma$
are completed by the punctures $p_1,\dots,p_n\in\bar{C}^\An$,
the completion of $C$. 
The induced tropicalisation map $C^\An\to\Gamma$ is the retraction map to the
skeleton.

\subsection{Deligne-Mumford compactification of $M^\trop_{g,n}$}

Let $2g-2+n>0$ and $\bar{\mathcal M}_{g,n}$ the moduli space
of $n$-pointed stable curves defined over $K^0$,
and $\bar{M}_{g,n}:=\bar{\mathcal M}_{g,n}\times_{\Spec K_0}\Spec K$.
 
For $S=\Spec K_0$ and an  $n$-pointed stable curve 
$(C\to S, s_1,\dots,s_n\colon S\to C)$ fix $e$ as in the previous
subsection, and  let $\Phi\colon C\to\mathbb{P}(V^*)\cong\mathbb{P}^N_S$
be the induced embedding, where
$V=H^0(C,\omega_{C/S}(s_1+\dots+s_n)^{\otimes e})$ 
and $N=e(2g-2+n)-g$. 

By \cite[Lemma 7.4]{Gubler1998}, the line bundle
$\mathscr{L}:=\omega_{C/S}(s_1+\dots+s_n)^{\otimes e}$ 
has a canonical formal metric $\norm{}_{\mathscr{L}}$.
This formal  metric can be applied to the projective embedding as follows.
Namely, take a basis $\sigma_0,\dots,\sigma_N$
of $V^*$ such that $\sigma_i$ restricted to the generic fibre $C_0$
of $C$ has a pole in
the $i$-th puncture, and let $L:=\mathscr{L}|_{C_0}$. 
The restricted 
 projective embedding $\Phi_L\colon C_0\to\mathbb{P}^N$
is defined on the chart 
$U_i=\mathset{x\in C_0\mid \sigma_i(x)\neq 0}$ as 
\begin{align}
x\mapsto \left(\frac{\sigma_1(x)}{\sigma_i(x)}:
\dots:1:\dots:\frac{\sigma_N(x)}{\sigma_i(x)}\right), \label{plurican}
\end{align}
By definition of the canonical formal metric, the restriction
$\norm{}$
of $\norm{}_{\mathscr{L}}$ to
$L$ satisfies
$-\log\norm{\sigma(x)}=v(\sigma(x))$
for $\sigma\in H^0(C_0,L)$.
Hence, on $U_i$,
$$
\Val(\Phi_L(x))=\left(v(\sigma_i(x))-v(\sigma_1(x)),\dots,v(\sigma_i(x))-v(\sigma_N(x))\right)\in\mathbb{TA}^N,
$$
and these maps glue to a map
$\Val\circ \Phi_L\colon C_0\to\mathbb{TP}^N$.
Note that we have used the property (\ref{can=taut})
and that the formal canonical metric
is the restriction of the formal canonical metric on $O(1)$.

Define $\trop(C):=\Val(\Phi_L(C))^\closure\subseteq{TP}^N$.
For $[C,p_1,\dots,p_n]\in\bar{M}_{g,n}$, we will label in
 the tropical curve $\Gamma=\trop(C)$  
the punctures $\Val(\Phi_{L_{C}}(p_i))$
by their original names $p_i$. This is reminiscent of the
fact that the skeleton 
 of a Berkovich-analytic space $X$
lies inside the space $X$.

\begin{rem}
If $(C,p_1,\dots,p_n)$ is an $n$-pointed projective line
defined over a finite extension $k$ of $\mathbb{Q}_p$, then
one can check that 
$\trop(C\times_k K)$ is isometric to the subtree
$\mathscr{T}^*(\mathset{p_1,\dots,p_n})$
of the Bruhat-Tits tree $\mathscr{T}_k$ of $\PGL_2(k)$, defined in \cite{Kato2005}.
That tree is the smallest subtree of $\mathscr{T}_k$ containing all geodesic
lines between the points $\mathset{p_1,\dots,p_n}\subseteq\mathbb{P}^1(k)$.
\end{rem}


We  propose a modified definition of
tropical curve in order to have a compact moduli space containing
the tropicalisations of $p$-adic curves. 

\begin{dfn}
Let 
$I:=I_1\coprod\dots\coprod I_n$ 
be a finite disjoint union of copies of the closed
or half-open unit
interval $[0,1]$ or $[0,1)$. Let the boundary
$\partial I=(I_1\setminus I_1)^\circ\coprod\dots\coprod (I_n\cup I_n^\circ)$
be partitioned into the disjoint union of non-empty sets 
$P_1\cup\dots \cup P_r$. 
The  topological quotient space
$\Gamma=I/\sim$ obtained by
the equivalence relation $\sim$ which
identifies the points in each $P_i$ is called a 
{\em semi-graph}. The equivalence classes 
of the $P_i$ are called {\em vertices}, and  the open intervals
$I_i^\circ$ {\em edges} of $\Gamma$. An edge $e$  is called {\em
bounded}, if its closure is the image of a closed interval under the canonical
projection $I\to\Gamma$. An edge which is not bounded is 
called a {\em puncture}.
\end{dfn}

\begin{dfn} \label{def-tropcurve}
A {\em tropical curve}  is a finite 
connected semi-graph $\Gamma$  whose vertices $v$
are labeled with natural numbers $g_v$, and whose bounded edges
are assigned lengths in $\mathbb{R}_{>0}\cup \mathset{\infty}$.  
A  tropical curve is {\em stable} if
 the vertices with label $0$ have at least three edges or
punctures attached to them,
and those with label $1$ have at least one edge or puncture emanating. 
A tropical curve is {\em smooth}, if it is stable and all bounded
edges have finite length.
The subsets $\Gamma^0, \Gamma^1,\Gamma^1_0,\Gamma^1_\infty$
of $\Gamma$ consist of the vertices, edges, bounded edges and punctures. 
\end{dfn}

The abstract tropical curves from \cite{GathmannMarkwigWDVV} are
smooth tropical curves according to Definition \ref{def-tropcurve}.

\begin{dfn}
The {\em arithmetic genus} of a  tropical curve $\Gamma$
is the number
$$
g(\Gamma):=1+\sum\limits_{v\in\Gamma^0}(g_v-1)+\#\Gamma^1_0.
$$
The {\em moduli space} of $n$-pointed stable tropical curves
of arithmetic genus $g$ is
$$
\bar{M}_{g,n}^\trop:=
\mathset{\Gamma\mid \text{$g(\Gamma)=g$ and $\#\Gamma^1_\infty=n$}}/\sim,
$$
where $\Gamma$ is a stable tropical curve, and
$\Gamma\sim\Gamma'$ means that there is a homeomorphism $\Gamma\to\Gamma'$
which takes vertices to vertices, bounded edges to bounded edges, 
punctures to punctures and respects the labellings.
\end{dfn}

\smallskip
We will in the following often suppress the labelling of vertices. However, it should
be born in mind that the interpretation is meant to be
that of  the labelling by
``genus of component of special fibre in stable model of a $p$-adic curve''.
The following theorem makes apparent the
  meaningfulness
of this interpretation. Later,  Corollary \ref{localgenus}
could be used to recover the suppressed labelling 
of vertices in finite covers of tropical curves.

\begin{thm} \label{barMgntrop}
The moduli space $\bar{M}_{g,n}^\trop$ is a 
compact polyhedral complex of pure dimension $3g-3+n$, and the
map
$$
\trop\colon\bar{M}_{g,n}\to\bar{M}_{g,n}^\trop,\;
[C,p_1,\dots,p_n]\mapsto[\trop(C),p_1,\dots,p_n]
$$
is well-defined and has dense image in $\bar{M}_{g,n}^\trop$.
\end{thm}

\begin{proof}
The space $\bar{M}_{g,n}^\trop$ is clearly a compactification
of the space $M_{g,n}^\trop$ of $n$-pointed smooth tropical curves 
by adding extra cells. That the space $M_{g,n}$ is a polyhedral complex
of the named pure dimension 
was seen in \cite[Ex.\ 2.13]{GathmannMarkwigWDVV} for $g=0$.
The general case can be proved in a similar way.

Let $\mathcal{X}\to\Spec{K^0}$ be a stable $n$-pointed curve.
The canonical formal
metric $\norm{}_{\mathscr{L}_{\mathcal X}}$ on the line bundle $\mathscr{L}$
from the beginning of this subsection
 restricts to
a metric $\norm{}$ on the line bundle
$L=\mathscr{L}|_X$ on the generic fibre $X$.
This  endows the stable reduction graph $\Gamma$ of $X$ with a metric. 
Namely, $(L,\norm{})$ induces a commuting
diagram
$$
\xymatrix{
X\ar[r]^{\Phi_{L}}\ar[d]&\mathbb{P}^N\ar[d]^{\Val}\\
\Gamma\ar[r]&\mathbb{TP}^N
}
$$
where $\Gamma$ is the combinatorial graph underlying the
closure of $\Val(\Phi_{L}(X))$ in $\mathbb{TP}^N$.

If $X$ is smooth, then  $\Gamma$ 
is given the structure of a tropical curve, as follows easily
with \cite[Thm.\ 2.1.1]{EKL2006}. 

Now assume that $X=X_0$ contains a double point  $x=x_0$.
We can realise $X_0$ as the fibre in $t=0$
of a  family $X_t$ of curves
flat over $K$ with constant reduction graph $\Gamma$,
and such that $X_t$ is smooth for $t\neq 0$.
Let $\mathcal{X}_t$ be the stable $K^0$-model of $X_t$,
and let $x_t\in X_t$ converge to $x_0$.
In other words, there is a section $s\colon T\to \mathfrak{X}$,
where $\mathfrak{X}\to T$ is the family $X_t$, 
such that $s(t)=x_t$. 
Let $\mathcal{L}:=\omega_{\mathcal{X}/T}(t_1,\dots,t_n)^{\otimes e}$, 
where $t_1,\dots,t_n\colon T\to\mathfrak{X}$ 
are the $n$ punctures of $\mathcal{X}\to T$,
and let $L_t:=\mathcal{L}|_{X_t}$. 
It is well known that the reduction map
$\pi_t\colon X_t^\An\to\mathcal{X}_{t,\sigma}$ 
to the special fibre $\mathcal{X}_{t,\sigma}$ 
of $\mathcal{X}_t$ has the property that 
the pre-image of a double point  $\xi_t\in\mathcal{X}_{t,\sigma}$ 
is an open annulus $A_t$.
The edge $e$ in $\Gamma$ corresponding to $\xi_t$
corresponds for $t\neq 0$ to an edge $e_t$ in 
$\Gamma_t:=\Val(\Phi_{L_{t}}(X_t))$ and
has length equal to the thickness of the annulus $A_t$.
Now assume that $A_t$ contains $x_t$. Because $x_t$ is not a puncture,
it defines a point $v\in e$ corresponding to a point
 $v_t\in e_t$ for $t\neq 0$.
Let $S_t\subseteq A_t$ be a circle containing $x_t$ such that
$\pi_t(S_t)=v_t$. We can realise the family $\mathfrak{X}$ in such
a way that $S_t$ shrinks for $t\to 0$ to the double point $x_t$,
while  the outer boundary circle of $A_t$ remains of constant radius $1$.
In this case, the length of $e_t$ increases unboundedly to $\infty$,
the length of $e_0$. By construction, the latter is the edge
corresponding to $e$
in the limit $\Gamma_0$ of the family of tropical curves $\Gamma_t$.
Hence, $\Gamma_0$ is a stable tropical curve.

The  arithmetic genus of 
the tropicalisation of $x\in\bar{M}_{g,n}$ is well-known
to
equal the arithmetic genus of any model
 of $x$ over $\Spec K^0$.
\end{proof}

The space 
$\bar{M}^\trop_{g,n}$
is also called the {\em Deligne-Mumford compactification of $M^\trop_{g,n}$}.
Important for us is 
a stratification of ${M}_{g,n}$
according to Betti number.

\begin{lem} \label{betti}
Let $C$ be an irreducible, non-singular projective curve of genus $g$.
Then for any constant sheaf $A$ on $C$, it holds true that
$H^1(C,A)\cong A^b$, where $b$ is the first Betti number of $\trop(C)$.
\end{lem}

\begin{proof}
\cite[Prop.\ 7.4.3]{FvP2004}.
\end{proof}

Lemma \ref{betti} implies that 
$$
\bar{M}_{g,n}=\coprod\limits_{b=0}^g\bar{M}_{g,n}^b,
$$
where $\bar{M}_{g,n}^b$ consists of the isomorphism classes
of punctured curves $(C,p_1,\dots,p_n)$ with $h^1(C)=b$.

\begin{cor}
The tropicalisation map
$
\trop\colon\bar{M}_{g,n}\to\bar{M}^\trop_{g,n}
$
satisfies 
$$
\trop(\bar{M}_{g,n}^b)\subseteq\bar{M}_{b,n}^\trop,
$$
where the image is dense.
\end{cor}

\begin{proof}
This is immediate from Proposition \ref{barMgntrop} and Lemma \ref{betti}.
Alternatively, a more careful  proof Proposition \ref{barMgntrop}
would yield the same result.
\end{proof}

It is known that $\bar{M}_{g,n}$ is a projective variety  
\cite[Thm.\ 6.1]{Knudsen1983}. It would be desirable to construct
a tropical ample line bundle on $\bar{M}_{g,n}^\trop$ in a 
similar way as can be done for its classical counterpart.

\subsection{Stable maps of degree $d$ to $\mathbb{P}^1$} \label{tropmor}

Let $f\colon\bar{C}\to\mathbb{P}^1$ be a stable map of degree
$d$ of $p$-adic 
curves which is defined over $K$.
We assume for $C:=\bar{C}\setminus\mathset{p_1,\dots,p_n}$
 that $[C,p_1,\dots,p_n]\in\bar{M}_{g,n}$. Let 
$B_f\subseteq \mathbb{P}^1(K)$
be the set of branch points of $f$. Then there exist $K^0$-models
$\mathcal{C}$ and $\mathcal{C}'$ of
$C$ and $C':=\mathbb{P}^1\setminus B_f$ such that $f$
extends to a map $F\colon \mathcal{C}\to\mathcal{C}'$.
This yields 
 an induced
map between tropical curves
$$
\trop(f)\colon\Gamma:=\trop(C\setminus f^{-1}(B_f))
\to\Gamma':=\trop(\mathbb{P}^1\setminus B_f)
$$
called the {\em tropicalisation} of $f$.



\section{A tropical Riemann-Hurwitz formula}

Let $\Gamma$ be a tropical curve. For reasons of convenience,
we assume $\Gamma$ to be smooth in the remainder of this article. 
It is often convenient for us to
deal with the tropical analogon of {\em complete} or {\em projective}
curve. For this reason, we adjoin to $\Gamma$ an additional vertex
at the end of each unbounded edge of $\Gamma$. This yields the
{\em completion} $\bar\Gamma$ of $\Gamma$. However, 
for us it is important to keep track of the set of ends
$\bar\Gamma\setminus\Gamma$ as they can be interpreted as $K$-rational
points on non-Archimedean curves for a suitable valued field $K$.

The {\em degree} $\deg (v)$ of a vertex $v\in\Gamma^0$ is defined as
$$
\deg(v):=\#\mathset{e\models v},
$$
where $e\models v$ means an edge $e$ emanating from a vertex $v$.
A tropical curve will be called  {\em binary}, if the degree of each
vertex is three. 
A {\em binary tree} is a binary tropical curve whose first Betti number
vanishes.

\begin{dfn}
Let $\Gamma$ be a tropical curve.  Then a {\em divisor} on $\Gamma$
is a formal finite sum
$$
D=\sum\limits_{v\in\bar{\Gamma}^0}D(v)[v]
$$
with $D(v)\in\mathbb{Z}$. The divisor  is expressed as $D=D_0+D_\infty$, 
where
$$
D_\infty:=\sum\limits_{v\in\bar{\Gamma}^0\setminus\Gamma^0}D(v)[v]
$$
is called the {\em part at infinity} of $D$.
The divisor
$$
K_\Gamma:=\sum\limits_{v\in\bar\Gamma^0}(\deg(v)-2)[v]
$$
is called the {\em canonical divisor} of $\Gamma$.
\end{dfn}

If $\Gamma$ is the tropicalisation of a smooth $p$-adic curve $X$,
then
the canonical divisor
$$
K_\Gamma=\sum\limits_{v\in\Gamma^0}(\deg(v)-2)[v]
-\sum\limits_{e\in\Gamma^1_\infty}[e]
$$
is the sum of a divisor on $X$ and a sum of vertices of $\Gamma$.
However, the part consisting of points of the completion $\bar{X}$
of
$X$ is different from
the canonical divisor of  $\bar{X}$. In fact, it
is merely minus the sum of  the punctures on $\bar{X}$ defined by $X$.
On the other hand, any divisor on $\Gamma$ can be interpreted
as an element of the free abelian group on the points of the
Berkovich analytification $\bar{X}^\An$ of the $p$-adic curve $\bar{X}$.

\begin{rem}
Clearly, for a tropical curve $\Gamma$, it holds true that
$$
\deg K_\Gamma=2b_1(\Gamma)-2,
$$
where $b_1(\Gamma)=b_1(\bar{\Gamma})=\#\bar\Gamma^1-\#\bar\Gamma^0+1$
is the first Betti number of $\Gamma$.
\end{rem}

A {\em morphism} $\phi\colon\Gamma\to\Gamma'$ of tropical
curves will mean for us in this article a continuous map between the underlying
topological spaces which takes edges to edges and vertices to vertices.

\begin{dfn}
A {\em weighted} tropical curve is a tropical curve $\Gamma$ together
with a map $w\colon\Gamma^1\to\mathbb{N}\setminus\mathset{0}$, 
called {\em the weights}. 
A morphism of tropical curves $\phi\colon\Gamma\to\Gamma'$ 
is said to be {\em $w$-harmonic}, if the weights
$w$ on $\Gamma$ are
such that for all $v\in\bar\Gamma^0$ the quantity
$$
m_{\phi,w}(v):=\sum\limits_{e\in\phi^{-1}(e'),\,e\models v}w(e)
$$
is independent of the choice of edge $e'\models \phi(v)\in\Gamma'^1$.
If $w$ is constant and equal to $1$, then a $w$-harmonic morphism will
be called {\em harmonic}. The quantity $m_{\phi,w}$ is called
the {\em multiplicity} of $v$ in $\phi$ with respect to $w$.
\end{dfn}

We will often suppress the weights $w$ in the notation of multiplicity,
when the weights are clear from the context.

\begin{lem} \label{degbound}
The following inequality holds true:
$$
\deg(v)\leq m_\phi(v)\deg(\phi(v))=:\deg_w(v),
$$ 
and is in general not an equality.
\end{lem}

\begin{proof}
The statement follows immediately from the obvious
equality
$$
\sum\limits_{e\models v}w(e)=m_\phi(v)\deg(\phi(v))
$$
and the fact that the values of the weights are at least one.
\end{proof}

The {\em degree} of a $w$-harmonic morphism generalises the
corresponding notion for graph morphisms from \cite{BN2007}:

\begin{dfn}
The {\em degree} of a $w$-harmonic morphism $\phi\colon\Gamma\to\Gamma'$
of tropical curves is defined as
$$
\deg\phi:=\sum\limits_{e\in\phi^{-1}(e')}w(e),
$$
where $e'\in\Gamma'^1$ is arbitrary.
\end{dfn}

\begin{rem}
Clearly, a morphism of degree $d$ is surjective.
We must, however, prove that the degree of $\phi$ is well-defined.
\end{rem}

\begin{proof}
The corresponding modification of the proof of \cite[Lemma 2.4]{BN2007} 
proves the statement.
\end{proof}

When we speak of a morphism of degree $d$, we usually mean 
that the weights of $\Gamma$ are fixed. There is some freedom
of choice due to the deck group of $\phi$, i.e.\ the
 automorphisms of $\Gamma$ leaving the fibres of $\phi$ invariant.

\begin{lem}
For any vertex $v'\in\Gamma'$, we have
$$
\deg(\phi)=\sum\limits_{v\in\phi^{-1}(v')}m_\phi(v).
$$
\end{lem}

\begin{proof}
The proof is as simple as that  of \cite[Lemma 2.6]{BN2007}. 
\end{proof}

\begin{dfn}
Let $\phi\colon \Gamma\to \Gamma'$ 
be a morphism of tropical curves of degree $d$.
Then
\begin{align*}
\phi_*\colon&\Div\Gamma\to\Div\Gamma',
\;D\mapsto\sum\limits_{v\in\bar\Gamma^0}D(v)[\phi(v)]\\
\phi^*\colon&\Div\Gamma'\to\Div\Gamma,
\;D'\mapsto\sum\limits_{v'\in\bar\Gamma'^0}\;
\sum\limits_{v\in\phi^{-1}(v')}m_\phi(v)D(v')[v]
\end{align*}
are the {\em push-forward} and the {\em pull-back} homomorphisms,
respectively.
\end{dfn}

Let $\phi\colon\Gamma\to\Gamma'$ be a morphism of tropical curves of degree 
$d$.
The divisor $R_\phi:=K_\Gamma-\phi^*K_{\Gamma'}$ is called the 
{\em ramification divisor} of $\phi$.
The {\em branch divisor} of $\phi$ is $\branch(\phi):=\phi_*R_\phi$.
A {\em tip} of a tropical curve $\Gamma$ is a vertex of degree one.
The set of tips of $\Gamma$ will be denoted by $\Gamma^0_\infty$.
A vertex $v\in\bar{\Gamma}^0$ such that $R_\phi(v)=0$ will be called 
{\em unramified}, just like in the classical case.

\begin{thm}[Tropical Riemann-Hurwitz formula] \label{RHtrop}
Let $\phi\colon\Gamma\to\Gamma'$ be a morphism of tropical curves of degree 
$d$.
Then 
it holds true that
\begin{align}
2b_1(\Gamma)-2&=(2b_1(\Gamma')-2)\cdot d+\deg R_\phi \label{RH}
\\
(R_\phi)_0
&=
\sum\limits_{v\in\Gamma^0}
  \left(2m_\phi(v)-2-(\deg_w(v)-\deg(v))\right)[v]. \label{ram_f}
\\
(R_\phi)_\infty&=\sum\limits_{v\in\bar\Gamma_\infty^0}(m_\phi(v)-1)[v]
\label{ram_infty}
\\
&b_1(\Gamma')\le b_1(\Gamma). \label{bettismaller}
\end{align}
In particular, the degree of $R_\phi$ is non-negative and even.
\end{thm}

\begin{proof}
Formula (\ref{RH}) follows easily from the formula
$$
\deg\phi^*D'=\deg\phi\cdot\deg D',
$$
whose proof is an adaptation of 
the unweighted version \cite[Lemma 2.13]{BN2007}.
This implies that $\deg R_\phi$ is of even degree. 
The formulae (\ref{ram_f}) and (\ref{ram_infty})
are straightforward.

Let for any graph $G$ denote $C(G)$ the chain complex associated to
$G$, after possibly assigning orientations to the edges. Then $\phi$
induces a morphism of chain complexes $C(\bar\Gamma)\to C(\bar\Gamma')$
which, by assumption, admits a well-defined
section $\sigma\colon C(\bar\Gamma')\to C(\bar\Gamma)$
satisfying
$$
[e']\mapsto \sum\limits_{e\in\phi^{-1}(e')}w(e)[e].
$$
The existence of the section  $\sigma$ implies
that there is an injective homomorphism
$
H_1(\Gamma',\mathbb{Z})\to H_1(\Gamma,\mathbb{Z})
$
from which (\ref{bettismaller}) follows.
\end{proof}

Let $f\colon C\to C'$ be a morphism of smooth 
projective $p$-adic curves
of degree $d$, and $\phi\colon\Gamma\to\Gamma'$ its
tropicalisation according to Section \ref{tropmor}. Assume that
$\phi$ is a morphism of tropical curves in our sense.
Let $v\in\Gamma^0$ and $v':=\phi(v)$. The {\em star} of $v$
denotes the tropical curve consisting of $v$
labelled with $g_v$ and the edges
in $\Gamma$ adjacent to $v$. The restricted
morphism
$\phi_v\colon\Gamma_v\to\Gamma'_{v'}$ is the tropicalisation
of a morphism $f_v\colon C_v\to C'_{v'}$ constructed as follows.
Via tropicalisation, the star $\Gamma_v$ corresponds to an affinoid
subset $U_v$ of $C$. It can be completed with affinoid disks to a projective
curve $C_v$ of genus $g_v$ \cite[Thm.\ 7.5.16]{FvP2004},
and one obtains by continuation of $f|_{U_v}\colon U_v\to U'_{v'}$
 a morphism $f_v\colon C_v\to C'_{v'}$
of degree $m_v\le d$. The branch locus $B_{f_v}$ is contained
in the set of punctures $\Gamma_{v'}^1$. Let $e'\models v'$. 
Then we assign to $e\in f_v^{-1}(e')$ as weight $w(e)$ the ramification degree
of the corresponding point in the cover $f_v$.
Note that $w(e)$ does not depend on the vertex $v$ adjacent to $e$.
Hence, we obtain weights $w$ on $\Gamma$ and call these the
{\em weights induced by $f$}. 

\begin{lem}
The tropical morphism as above is $w$-harmonic for
the  weights $w$ on $\Gamma$ induced by the morphism $f$.
\end{lem}

\begin{proof}
This follows from the local equality
$$
\sum\limits_{e\in f_v^{-1}(e')}w(e)=m_v
$$
for each vertex $v\in\Gamma$.
\end{proof}

Whenever $\Gamma\to\Gamma'$
is the tropicalisation of a finite map of smooth $p$-adic curves,
we assume that the weights of $\Gamma$ are induced by $f$.

\begin{thm}[RHM-criterion] \label{MRH}
Assume that  $\phi\colon\Gamma\to\Gamma'$ is
a morphism of tropical curves  which is 
the tropicalisation of a 
map $f\colon X\to X'$ of smooth projective
$p$-adic curves of degree $d$, where $X'$ is a Mumford
curve. 
Then the following statements are equivalent:
\begin{enumerate}
\item $X$ is a Mumford curve. \label{XisMumf}
\item $R_\phi$ is effective.  \label{Rphi>0}
\item $R_\phi=R_f$, where $R_f$ is the ramification 
divisor of $f$.               \label{Rphi=Rf}
\end{enumerate}
\end{thm}

\begin{proof}
(\ref{Rphi=Rf}) $\Rightarrow$ (\ref{XisMumf}). This follows from
Theorem \ref{RHtrop} and the classical Riemann-Hurwitz theorem for curves.

(\ref{XisMumf}) $\Rightarrow$ (\ref{Rphi=Rf}).
Now assume that $X$ is a Mumford curve. First, consider the case that
$\Gamma$ has precisely one vertex $v$. 
Then $R_\phi(v)=0$, as otherwise $\deg R_\phi\neq \deg R_f$ yields the
contradiction $g(X)\neq b_1(\Gamma)$ (where the latter is $0$).
In the general case,  the $\Gamma_v$
 of a vertex $v\in\Gamma^0$
yields as above a morphism of tropical curves associated to the 
cover $f_v\colon X_v\to X'_{\phi(v)}$ of projective smooth
$p$-adic curves.
From the assumption, it follows that $X_v$ and $X'_{v'}$  
are Mumford curves.
By the first case, it holds true
that $R_\phi(v)=R_{\phi_v}(v)=0$. 
It follows that $R_\phi=R_f$, as claimed.

(\ref{Rphi=Rf}) $\Rightarrow$ (\ref{Rphi>0}). This is obvious.

(\ref{Rphi>0}) $\Rightarrow$ (\ref{Rphi=Rf}). Note that for $v\in\Gamma^0$
it holds true that
\begin{align}
2m_\phi(v)-2+2g(X_v)=\deg R_{f_v}&=\sum\limits_{x\in\bar\Gamma^0_{v,\infty}}(w(e)-1)
=\deg_w(v)-\deg(v)
\label{localRH}
\end{align} 
where the first equality holds true by the classical Riemann-Hurwitz theorem.
Hence, $R_\phi(v)> 0$ implies (by Theorem \ref{RHtrop}) 
$\deg R_{f_v}< 2m_\phi(v)-2$, from which it follows by (\ref{localRH}) 
that $g(X_v)< 0$.
This is a contradiction.
 So, if $R_\phi$ is effective, then necessarily
$R_\phi=R_f$.
\end{proof}

We remark that, due to Theorem \ref{MRH}, we cannot
define  a straightforward 
tropical version of the Fantechi-Pandharipande map
$\bar{M}_{g,0}(\mathbb{P}^1,d)\to\mathbb{P}^r$, where $r$ is determined by 
the classical Riemann-Hurwitz theorem. What we can get
for general tropical  curves is a map from 
$\bar{M}_{g,0}(\mathbb{P}^1,d)$ into the monoid of effective
divisors on  $(\mathbb{P}^1)^\An$. 

\begin{cor} \label{localgenus}
In the situation as in Theorem \ref{MRH}
it holds true that
$$
R_\phi(v)=-2g(C_v)=-2g_v,
$$
where $v\in\Gamma^0$ and $C_v$ is defined after the proof of 
Theorem \ref{RHtrop}.
\end{cor}

\begin{proof}
This follows easily from (\ref{localRH}).
\end{proof}

\begin{rem}
The result of Theorem \ref{MRH} is that the degree of the ramification
divisor $R_\phi$ depends on the relative position of the
branch points of the covering $f\colon X\to X'$. 
\end{rem}

\section{Tropical Hurwitz numbers}

In this section, we develop tropical methods for counting Mumford
curves covering the $p$-adic projective line depending
on the position of the branch points. This can be viewed as an 
extension
of \cite{Brad-MathZ251, Brad-manmat124}, where the focus
was on Galois covers
of Mumford curves.
 
\subsection{Binary branch trees}

After dealing with  comb-shaped targets, we prove that tropical Hurwitz numbers
above binary trees count solely Mumford curves.

\begin{dfn} \label{tropcomb}
By a {\em tropical comb} we mean a tropical curve $\Gamma$ of the form:
$$
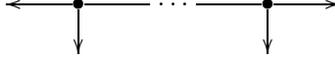
\xymatrix@R=15pt{
\ar@{<-}[r]&*\txt{$\bullet$}\ar[d]\ar@{-}[r]&\dots\ar@{-}[r]
&*\txt{$\bullet$}\ar[r]\ar[d]
&\\
&&&
}
$$
The horizontal geodesic line will be called the  {\em backbone},
the two horizontal ends the
 {\em handles}, and the other ends the {\em teeth} of $\Gamma$.
\end{dfn}

Let a  set 
$\mathset{\eta_0=(\eta_0^1,\dots,\eta_0^{\ell_0}),
\dots,\eta_{n+1}=(\eta_{n+1}^1,\dots,\eta_{n+1}^{\ell_{n+1}})}$ of partitions of a
natural number $d$, and   a tropical comb $\Gamma'$ with $n$ vertices
$v_1,\dots,v_n$, be given, where the number $n+2\ge 3$ is 
the number of branch points defined by the classical
Riemann-Hurwitz formula.
We will construct covers $\Gamma\to\Gamma'$ according to the following rules:

\begin{enumerate}
\item Start with $\ell_0$ ends above the left handle, and weight these
with the entries of $\eta_0$. These ends will be called {\em strands}.
\item Above $v_1$ create vertices  by joining and splitting or continuing 
strands above the next edge of the backbone.
Weight the new outgoing strands such that for each vertex $v$
 the sum of incoming weights 
equals the sum of outgoing weights. This number is $m_v$
\item From all vertices $v$ above $v_1$ let strands emanate above the tooth
$t_1$ attached to $v_1$ such that the total number of strands above
$t_1$ is $\ell_1$, their weights are the entries of $\eta_1$,
and the sum of the weights of the new 
strands emanating from $v$ equals again $m_v$. 
\item In (2) and (3) make sure that  all vertices $v$ above $v_1$
satisfy $\deg(v)=m_v+2$.
\item Repeat the procedure successively above the vertices $v_2,\dots,v_n$,
observing rule (4) each time.
\item Above the right handle make sure there are $\ell_{n+1}$ strands
weighted with the entries of $\eta_{n+1}$.
\end{enumerate}

\begin{dfn}
The number $H_d^g(\eta_0,\dots,\eta_{n+1})^\trop$ is defined to be the
number of isomorphism classes of covers $\Gamma\to\Gamma'$
satisfying the above rules $(1)$ to $(6)$ with $b_1(\Gamma)=g$.
It is called {\em tropical Hurwitz number}.
\end{dfn}

Notice, that
if we set $\eta_0=\eta$, $\eta_{n+1}=\nu$ and $\ell_i=d-1$
for $i=1,\dots,n$, then $H_d^g(\eta_0,\dots,\eta_{n+1})^\trop$
specialises to the double Hurwitz number $H_d^g(\eta,\nu)$
from \cite{CJM}. In fact, the authors of \cite{CJM} consider covers
of the tropical projective line, which can be obtained from our covers
by contracting all teeth. The rules (2) to (4) 
applied to the case of
simple ramification enforce the vertices of the upper tropical curve
to have precisely three adjoint strands above the backbone of the comb.

\smallskip
Let $H_d^g(\eta_0,\dots,\eta_{n+1})$ denote the Hurwitz
number of smooth connected covers of $\mathbb{P}^1$ of degree $d$ and
genus $g$ ramified above $n+2$ points with ramification profiles
$\eta_0,\dots,\eta_{n+1}$. If we fix the branch locus 
$B\subseteq\mathbb{P}^1(K)$, then we can  associate to
the lower curve the moduli point 
$x=[\mathbb{P}^1\setminus B,B]\in M_{0,n-1}$.
The quantity
$H_d^g(\eta_0,\dots,\eta_{n+1})^\Mumf(x)$ then
denotes the number of Mumford curves of genus $g$ covering 
any representative of $x$
 in degree $d$ with the assigned ramification
profiles. Taking into consideration Remark \ref{preciseedgelengths} 
below,
this restricted Hurwitz number depends only on the combinatorial
type 
of the tropicalisation of the lower punctured curve. 
The inequality 
$H_d^g(\eta_0,\dots,\eta_{n+1})^\Mumf(x)\le H_d^g(\eta_0,\dots,\eta_{n+1})$
trivially holds true.

\begin{rem} \label{preciseedgelengths}
When we speak in the following of certain topological types
of tropical curves, we require that all bounded edges have some minimal
length $\ge 0$, depending on whether or not the prime $p$ divides
the ramification index of an adjacent vertex. The reason is
that only edges of sufficient length can possible come from 
Mumford curves; in the contrary case, there fails to be a
discrete action on the 
corresponding Bruhat-Tits tree.
The precise minimal edge
lengths can be calculated with the methods from \cite{Brad-manmat124}.
\end{rem}

\begin{thm} \label{combtrop=mumf}
It holds true that
$$
H_d^g(\eta_0,\dots,\eta_{n+1})^\trop=H_d^g(\eta_0,\dots,\eta_{n+1})^\Mumf(x)
=H_d^g(\eta_0,\dots,\eta_{n+1})
$$
where $x\in M_{0,n-1}$ is in comb position.
\end{thm}

\begin{proof}
We need to prove that $H_d^g(\eta_0,\dots,\eta_{n+1})^\trop=
H_d^g(\eta_0,\dots,\eta_{n+1})$. Then the 
first equality is a consequence of Theorem \ref{MRH},
as for any cover $\phi\colon\Gamma\to\Gamma'$ according to the rules
(1) to (6), the ramification divisor
$R_\phi$ is effective by construction. Clearly, the branch locus of
any cover of $p$-adic curves with tropicalisation equal to $\phi$
is in comb position.

Let $\phi\colon\Gamma\to\Gamma'$ be a tropical cover counted by
$H_d^g(\eta_0,\dots,\eta_{n+1})$. It can be obtained by glueing
local covers $\phi_v\colon\Gamma_v\to\Gamma_{\phi(v)}'$ defined
after  the proof
of Theorem \ref{RHtrop}. The trees 
$\Gamma'_{\phi(v)}$ consist of the vertex $v_i=\phi(v)$ and three ends,
one of which is the tooth $t_i$.
For fixed $v_i'$  the cover 
$$
\phi_i\colon\Gamma_i:=\coprod\limits_{v\in\phi^{-1}(v_i')}\Gamma_v
\to\Gamma'_{v_i'}
$$
is of ramification type $(\gamma_{i}^-,\eta_i,\gamma_{i}^+)$,
and the glueing condition is $\gamma_{i}^+=\gamma_{i+1}^-$.
Of course, $\gamma_0^-=\eta_0$ and $\gamma_{n+1}^+=\eta_{n+1}$.
The covers $\phi_v$ correspond to  covers 
$f_v\colon\mathbb{P}^1\to \mathbb{P}^1$ of the 
same degrees and ramification types.
Hence, $H_d^g(\eta_0,\dots,\eta_{n+1})^\trop$ counts
all possible genus $g$ covers of the projective line obtainable by  
glueing genus $0$ covers of $\mathbb{P}^1$ according to the same
glueing conditions. By specialising the Degeneration Theorem 
(or join-cut recursion)
\cite[Thm.\ 3.15]{Li2002}, it follows that 
this method counts all the way up to $H_d^g(\eta_0,\dots,\eta_{n+1})$.
\end{proof}

As a corollary, we obtain that double tropical Hurwitz numbers
count Mumford curves:

\begin{cor} \label{comb=mumf}
It holds true that
$$
H_d^g(\eta,\nu)^{CJM}=H_d^g(\eta,\eta_1,\dots,\eta_n,\nu)^\Mumf(x),
$$
where $H_d^g(\eta,\nu)^{CJM}$ is the double Hurwitz number
from \cite{CJM}, $\ell(\eta_i)=d-1$, and
$x\in M_{0,n-1}$ is in comb position.
\end{cor}

\begin{proof}
Let $h\colon\Gamma\to\mathbb{R}\cup\mathset{\pm\infty}$ be a tropical
map to $\mathbb{P}^1$ of degree $\Delta$
in the sense of \cite{CJM}, where the multiset $\Delta$
is given as
$$
\Delta=\mathset{-\eta_1,\dots,-\eta_{\ell(\eta)},\nu_1,\dots,\nu_{\ell(\nu)}}.
$$
Let $v'\in\mathbb{R}$ be a branch point of $h$. We may assume that
 $h^{-1}(v')$ is a discrete subset of $\Gamma$.
By subdivision of edges, we may further assume that every pre-image
of $v'$ is a vertex. Now, attach to every branch point in $\mathbb{R}$
a half-line in order to  obtain a tropical comb.
Attach also to every vertex $v$ of  $\Gamma$ a set of $m(v)-1$ half-lines,
where $m(v)$ is the sum of the weights on edges coming into $v$ from the left.
The resulting map is an extension of $h$ to 
a morphism $\phi\colon\Gamma\to \Gamma'$ of degree $d$ 
simply ramified above the teeth, and ramified with profiles $\eta$ and
$\nu$ above $\pm\infty$. Clearly, $(R_\phi)_0=0$.
Hence,  
$H_d^g(\eta,\nu)^{CJM}$ counts covers of $p$-adic $\mathbb{P}^1$ 
of degree $d$ by curves $C$ such that $g(C)=b_1(\trop(C))$, i.e.\
Mumford curves. Hence, 
$H_d^g(\eta,\nu)^{CJM}= H_d^g(\eta,\eta_1,\dots,\eta_n,\nu)^\trop$
which, by Theorem \ref{combtrop=mumf} counts Mumford curves 
ramified over points in comb position.
\end{proof}

For $x=[C,p_1,\dots,p_n]\in M_{0,n-1}$, denote  
 $[\trop(C),p_1,\dots,p_n]$
simply as $\mathscr{T}(x)$.

\begin{thm} \label{bin=mumf}
Let $x\in M_{0,n-1}$ such that $\mathscr{T}(x)$ is binary.
Assume that $\eta_0,\dots,\eta_{n+1}$ are $n+2$ partitions of $d$.
Then 
$$
H_d^g(\eta_0,\dots,\eta_{n+1})^\Mumf(x)
=H_d^g(\eta_0,\dots,\eta_{n+1}).
$$
\end{thm}

\begin{proof}
The binary tree $T=\mathscr{T}(x)$ can be obtained by glueing
combs $T_i=\mathscr{T}(x_i)$ along an extra end $e_i$ for each 
glueing morphism. 
The consequence is a construction of 
tropical degree $d$ maps $f\colon \Gamma\to\Gamma'$ from genus $g$ curves via 
local pieces $f_i\colon\Gamma_i\to\Gamma_i'$ 
from genus $g_i$ curves 
and the end $e_i$ corresponds to a  branch point
of $f_i$.  Clearly, the genera sum up to at most $b_1(\Gamma)\le g$.
Since the vertices in $\Gamma_i$ above the origin of 
$e_i$ are unramified by Lemma \ref{comb=mumf} and Theorem \ref{MRH}, 
the same holds true for the vertices of $\Gamma$.
Putting together 
$b_1(\Gamma_i)=g_i$ and  Theorem \ref{MRH} again, it follows 
that  the Hurwitz number does indeed count Mumford curves above $x$.

Using again the classical join-cut recursion as in the proof of 
Theorem \ref{combtrop=mumf}, it follows that the tropical
count equals the classical count.
\end{proof}

An immediate consequence is:

\begin{cor} \label{trop-reorder}
The number 
$H_d^g(\eta_0,\dots,\eta_{n+1})^\trop$ does not depend on the
ordering of the partitions $\eta_i$.
\end{cor}

\subsection{Edge Contractions}

Let $x\in M_{0,n-1}$ and $T=\mathscr{T}(x)$ the associated tropical curve.
Assume that $T$ is a binary curve, and contract a bounded edge $e'\in T^1$
in order to obtain the tropical curve $\Gamma'$.
Let $x'\in M_{0,n}$ be a punctured projective line such that
$\mathscr{T}(x')=\Gamma'$.
Then $H_d^g(\eta,\nu)$ does not count anymore only Mumford curves 
covering $x'$. Namely, any tropical curve $\Gamma$ with $b_1(\Gamma)=g$
above $\Gamma'$ having the interior of a graph shaped 
like\footnote{This is called {\em wiener} in \cite{CJM}.}
$$
\xymatrix{
*\txt{$\bullet$}\ar@{-}@/^/[r]\ar@{-}@/_/[r]&*\txt{$\bullet$}
}
$$
 above an edge $e'$ will be contracted
to the tropicalisation of a curve $C$ having the property $b_1(\trop(C))<g$.
Hence, for the double Hurwitz numbers we obtain the result:
\begin{align}
H_d^g(\eta,\eta_1,\dots,\eta_n,\nu)^\Mumf(x')=H_d^g(\eta,\nu)^{CJM}
-W(e') \label{nowieners}
\end{align}
where $W(e')$ is the number of double Hurwitz covers of degree $d$ and genus
$g$ such that $\phi^{-1}(\bar{e}')$ contains a wiener. Here,
$\bar{e}'$ is the segment in $T$ consisting of $e'$ and its adjacent
vertices.

In the  case of general ramification types above combs, 
we must let $W(e')$ count {\em multi-wieners}
above $e'$ in order to obtain the correct count in (\ref{nowieners})
for Mumford curves. By a multi-wiener, we mean  a graph of the following kind:

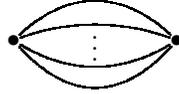
\vspace*{1mm}
$$
\xymatrix{
*\txt{$\bullet$}\ar@{-}@/^/[rr]\ar@{-}@/^15pt/[rr]
\ar@{-}@/_10pt/[rr]\ar@{-}@/_18pt/[rr]
&\vdots&*\txt{$\bullet$}
}
$$

\bigskip
This gives an algorithm for computing Hurwitz numbers for Mumford curves
covering $\mathbb{P}^1$ with given branch locus in $M_{0,n+2}$.

\subsubsection{An algorithm}
Let $\Gamma'=\mathscr{T}(x)$, where $x\in M_{0,n-1}$ is the branch locus.
Let $\gamma$ be a shortest path in $M_{0,n-1}^\trop$ to the nearest cell of 
maximal dimension (it consists of binary trees)
successively  passing  through cells of dimensions increasing by one.
At the other end of the path $\gamma$ lies a point $y\in M_{0,n-1}$
such that $T=\mathscr{T}(y)$ is a binary tree. 
The path induces a choice of edges
$E_\gamma=\mathset{e_1,\dots,e_{\ell(\gamma)}}$
in $T$ which get contracted in order to obtain $\Gamma'$.
Let $W(E_\gamma)$ 
be the number of Hurwitz covers $\phi\colon\Gamma\to T$ counted by 
$H_d^g(\eta_0,\dots,\eta_{n+1})^\trop$ such that 
$b_1(\phi^{-1}(\bar{E}_\gamma))>0$,
where $\bar{E}_\gamma$ is the closure of $E_\gamma$ in $T$.   

\begin{prop}
The quantity 
$W(E_\gamma)$ does not depend on the choice of path the $\gamma$.
It holds true that
$$
H_d^g(\eta_0,\dots,\eta_{n+1})^\Mumf(x)=H_d^g(\eta_0,\dots,\eta_{n+1})^\trop
-W(E_\gamma).
$$ 
\end{prop}

\begin{proof}
This follows  immediately from Theorem \ref{bin=mumf}
and the fact that the classical Hurwitz number is independent
of the positions of the branch points.
\end{proof}

\subsection{Cyclic covers: the Harbater-Mumford condition}

Assume that $f\colon C\to\mathbb{P}^1$ is a cyclic cover of degree $n$.
The methods we have developped allow to recover the well-known
{\em Harbater-Mumford} condition which are necessary for $C$ to be a
Mumford curve. Namely, if $C$ is a Mumford curve, then
the local tropical pieces of the cover are
themselves cyclic covers $\Gamma_v\to\Gamma_{v'}$
of degree $m_v$.
By  tropical Riemann-Hurwitz, the degree of the
vertex $v$ must be $m_v+2$. Hence, there can be ramification only
above two of the three ends, and each ramification degree
is $m_v$.
It follows that the corresponding cover of $p$-adic curves
$f_v\colon \mathbb{P}^1\to\mathbb{P}^1$ is cyclic of degree $m_v$,
and the two branch cycles $\sigma$ and $\tau$ are inverse to each other.

However, the Harbater-Mumford condition is not sufficient.
In fact, the previous section shows that the branch points
must form a tree whith  sufficiently many edges
in order for the upper curve to have the correct first Betti number.
In \cite{Brad-manmat124}, we have calculated the precise
branch trees (depending also on the prime number $p$) for which 
cyclic covers 
have Mumford curves on the top.

\begin{exa}
Let $E$ be an elliptic curve, and $f\colon E\to\mathbb{P}^1$ the
cover of degree $2$ which we may assume to be ramified in the points
$0,1,\lambda,\infty$
with $v(\lambda)=0$. It is known that $E$ is a Tate curve, if and only if
$v(\lambda-1)> 2\cdot v(2)$ (cf.\ \cite[Thm.\ 5]{Tate1974} for $p\neq 2$
and \cite[Ex.\ 3.8]{Brad-manmat124} for $p$ arbitrary). In fact,
after making the edge in the branch tree long enough, the upper tropical
curve $\trop(E\setminus f^{-1}\mathset{0,1,\lambda,\infty})$ 
obtains a wiener, the length of which can be viewed as the
tropical $j$-invariant of the tropical elliptic curve by the well-known 
formulae
 relating $j$, $\lambda$ and the Tate-parameter $q$
for which $E$ has the representation as Tate curve
$E\cong\mathbb{G}_m/q^\mathbb{Z}$.

In the case $v(j)\le 4\cdot v(2)$, there is still a
corresponding moduli point in $\bar{M}_{1,1}^\trop$, namely the
point representing the tree $\xymatrix{\ar@{-}[r]&*\txt{$\bullet$}}$.

The other extreme point in $\bar{M}_{1,1}$ corresponds to 
the rational curve with one self-intersection point. Its
tropicalisation has a loop of infinite length.
Hence, we have exhibited
$\bar{M}_{1,1}^\trop$ as a cw-complex isomorphic
to the segment 
$\xymatrix{
*\txt{$\bullet$}\ar@{-}[r]&*\txt{$\bullet$}
}$
and whose edge  contains densly the tropicalisations of 
Tate curves.
\end{exa}




\end{document}